\documentclass[11pt]{amsart}
\usepackage{amssymb,amsmath,amsthm,epsf,epsfig,bbm}

\usepackage{latexsym}
\usepackage{upref, eucal}
\usepackage[all]{xy}

\newcommand {\nc} {\newcommand}
\newcommand {\enm} {\ensuremath}

\nc {\bdm} {\begin{displaymath}}
\nc {\edm} {\end{displaymath}}

\newtheorem {theorem} {Theorem}[section]
\newtheorem {lemma} [theorem]{Lemma}

\newtheorem {remark}{Remark}[section]

\newtheorem {corollary} [theorem] {Corollary}
\newtheorem {proposition}[theorem]{Proposition}

\numberwithin{equation}{section}

\nc{\J}{\enm{\mathcal{J} }}
\nc {\Z} {\enm{\mathbb{Z}}}

\nc {\stk} {\stackrel}
\newcommand{\imp}{\Rightarrow}
\newcommand{\map}{\rightarrow}

\newcommand{\beqar}{\begin{eqnarray*}}
\newcommand{\eeqar}{\end{eqnarray*}}

\newcommand{\Pn}[2] {\ensuremath{ {\mathbb{P}}^{#1}_{#2}}}
\nc{\Quot}[3]{\enm{ {\mathfrak{Quot}_{ {#1}/{#2}/{#3}}}}}
\nc{\Hilb}[2]{\enm{ {\mathfrak{Hilb}_{ {#1}/{#2}}}}}

\newcommand{\bb}[1]{\mathbb{#1}}
\newcommand{\mcal}[1]{\mathcal{#1}}

\nc {\Coh}[4] {\ensuremath{H^{#1}(\Pn{#2}{},{#3}({#4}))}}
\nc {\Ch}[3] {\enm{H^{#1}(X_t,{#2}_t({#3}))}}
\nc {\Qphi}[4]{\enm{ {\mathfrak{Quot}^{~#4}_{ {#1}/{#2}/{#3}}}}}
\nc {\Gra}[4]{\enm{ {\mathfrak{Grass}_{#2}({#3},{#4})}}}

\title{ \bf Invariant Forms for Correspondences of Curves}
\author {Arnab Saha}
\date {}

\begin {document}
\maketitle

\section{\bf Introduction}

Consider the correspondences of curves defined over either an algebraically
closed field or a number field. In this article, we try to classify the 
invariant differential forms which a correspondence admits. To state things
precisely, let us recall that a correspondence in any category $\mcal{D}$ is a 
tuple $\mathbb{X} = (Y,X,\sigma_1, \sigma_2)$ where $X$ and $Y$ are objects in
$\mcal{D}$ and $\sigma_1,~\sigma_2: X \map Y$ are morphisms. In our case, 
$\mcal{D}$ is the category of smooth algebraic curves defined over $k$, where 
$k$ is either an algebraically closed curve of any characteristic or $k$ is a 
number field. Let $\Omega_{K(Y)}$ denote the sheaf of rational $1$-forms over 
the curve $Y$ with function field $K(Y)$. $\Omega_{K(Y)}^{\otimes \nu} $ will 
denote its higher tensor powers of $\Omega_{K(Y)}$ for any $\nu \in \Z$. 
We would say that a form $\omega \in \Omega_{K(Y)}^{ \otimes \nu},~\omega 
\ne 0$ is {\it invariant} of weight $\nu$ in $\mathbb{X}$ if 
$\sigma^*_1 \omega = \sigma^*_2 \omega$. A form $\omega$ will be called 
{\it sem-invariant} if $\sigma^*_1 \omega = \lambda \sigma^*_2 \omega$ for 
some $\lambda \in k^\times$.

Now, given a correspondence $\mathbb{X}$, we can associate the group 
$\mathcal{G}_\mathbb{X}$ defined as
$$\mathcal{G}_{\bb{X}} := \{\omega~|~\omega \ne 0,~ \sigma_1^* \omega = \lambda 
\sigma_2^* \omega,~\lambda \in k^\times \}/\sim$$ 
\noindent
where the equivalence
relation $\sim$ is $\omega \sim \omega'$ if and only if $\frac{\omega}{\omega'}
$ is a constant function on $Y$. If the degrees of the morphisms $\sigma_1$ and
$\sigma_2$ are unequal, then Proposition \ref{cyclic} implies that
$\mcal{G}_\mathbb{X}$ is in fact free group of rank $\leq 1$.
The first question that we would like to address is that, if the group 
$\mcal{G}_\mathbb{X}$ is non-trivial, then we would try to classify the 
differential forms $[\omega] \in \mcal{G}_\mathbb{X}$. Given a form $\omega$
on the curve $Y$, we define the {\it Conductor} of $\omega$, denoted by 
$\mcal{C}_{Y,\omega}$, to be the cardinality of the support of the divisor 
$div~\omega$ on $Y$, that is $\mathcal{C}_{Y,\omega} = \#(\mbox{support}~
(div~\omega))$. For example, if $Y = \mathbb{P}^1$ and $\omega = \frac{dt}{t}$
then $div~\omega = -0 -\infty$ and hence $\mcal{C}_{\mathbb{P}^1,\omega} = 2$.
Our first result on general correspondences of curves $\mathbb{X}$ is 
a bound on the conductor of a semi-invariant form $\omega$ in the case when
the degrees of the morphisms $\sigma_1$ and $\sigma_2$ are unequal. The 
result says that $\mcal{C}_{Y,\omega}$ can be uniformly bounded in terms 
of the genus of our projective, smooth, curves $X$ and $Y$ and the degrees of 
the morphisms $\sigma_1$ and $\sigma_2$. It is noteworthy to mention that this 
bound is independent of the characteristic of the ground field $k$.

\begin {theorem}
\label{bound}
Let $\mathbb{X}= (Y, X,\sigma_1, \sigma_2)$ be a correspondence defined 
over an algebraically closed field $k$, where $X$ and $Y$ are projective 
smooth curves with $deg~\sigma_1 > deg~\sigma_2$,
$\sigma_1$ and $\sigma_2$ are tamely ramified and let $\omega$ be a 
semi-invariant form. Then the conductor $\mathcal{C}_{Y,\omega}$ is bounded by 
a number entirely determined by the genus of $X$ and $Y$ and the degrees of the
 morphisms $\sigma_1$ and $\sigma_2$. 
\begin{equation}
\mathcal{C}_{Y,\omega} \leq \frac{1}{d_1-d_2}(3(2g_X-2) -(2d_1+d_2)(2g_Y-2)),
\end {equation}
where $deg~\sigma_i = d_i,$ for $i = 1,2$ and $g_X$ and $g_Y$ denote the genus 
of the curves $X$ and $Y$ respectively.
\label {thm1} 
\end {theorem}

The above theorem is only true in the case when the degree of the maps 
$\sigma_1$ and $\sigma_2$ are unequal. If they are equal, then the above
result is false, that is, there can not be any bound on 
$\mcal{C}_{Y,\omega}$. We will exhibit an example on Hecke 
Correspondences which shows that the conductor $\mcal{C}_{Y,\omega}$
grows linearly with the characteristic of the ground field.
Consider the Hecke correspondences \cite{Diamond},
$$ \mathbb{X} = (X_1(1), X_1(1,l), \sigma_1,\sigma_2) \mbox{  over } 
\mathbb{Z}[1/l], $$ where $X_1(1)$ and $X_1(l)$ are modular curves of level
$1$ and $l$ respectively and $\sigma_1$ and $\sigma_2$ are the 
degeneracy maps.  Then for each prime $p$, one can consider the correspondence
$\bb{X}_p = (X_1(1)_p, X_1(1,l)_p, \bar{\sigma}_{1,p},\bar{\sigma}_{2,p})$ 
obtained by base changing $\bb{X}$ to the algebraic closure of the residue 
field at the prime $p$, $ p \nmid l$. Then for each $p$, there exists an 
invariant form $\omega_{p}$ such that 
$(\mbox{support }div~(\omega_{p})) = \mbox{super-singular elliptic curves}$ 
\cite{B2}. Hence $$\mathcal{C}_{X_1(1)_p,\omega_p} = \#\{\mbox{distinct roots 
of the Hasse polynomial}\} = \left[\frac{p}{12}\right] + 2$$ because of the 
irreducibility of the Hasse polynomial \cite {Katz} and this shows that  
$\mathcal{C}_{X_1(N)_p, \omega_p}$ is unbounded. Therefore Theorem 
\ref{bound} shows us that having unequal degrees of the maps makes the 
conductor bounded uniformly.

We will call a non-zero form $\omega$ on $Y$ {\it primitive} if its class 
$[\omega]$ belongs to$\mcal{G}_\bb{X}$ and also generates $\mcal{G}_\bb{X}$. 
From now on, we will consider correspondences of the projective line $\bb{P}^1$
that is we have $X=Y = \bb{P}^1$. We define a pair of tuples of morphisms
$(\sigma_1, \sigma_2)$ and $(\sigma_1',\sigma_2')$ to be {\it conjugates}
if there exists an automorphism $\varphi$ of $\bb{P}^1$ satisfying
$\sigma_1' = \varphi \circ \sigma_1 \circ \varphi^{-1}$ and $\sigma_2' = 
\varphi \circ \sigma_2 \circ \varphi^{-1}$. We define a form $\omega$ on
$\bb{P}^1$ to be a {\it flat} form if it is either $\omega = \frac{dt}{(t-a)}$ 
or $\omega = \frac{(dt)^2}{(t-a)(t-b)}$ for some $a$ and $b$ and $a \ne b$.
Note that the form $\frac{dt}{(t-a)}$ has weight $1$ where as the form
$\frac{(dt)^2}{(t-a)(t-b)}$ has weight $2$.
Let us consider a correspondence $\bb{X} = (\bb{P}^1, \bb{P}^1,\sigma_1,
\sigma_2)$ where both $\sigma_1$ and $\sigma_2$ are tamely ramified. Our 
next result classifies the types of semi-invariant forms that the above
correspondence $\bb{X}$ can admit. However, we had to assume a technical
condition of sufficient seperateness between the degrees of the morphisms.
\begin {theorem}
\label{mainthm2}
Let $\bb{X}$ be as above and 
$(\sigma_1,\sigma_2) \sim (\sigma_1',\sigma_2')$ where $\sigma_1'$ and
$\sigma_2'$ are completely ramified at a point. Also further assume that
$deg~\sigma_1 \geq 14~ deg~\sigma_2$. If $\mathcal{G}_{\bb{X}}$ is 
non-trivial and $\omega$ is its primitive, then $\omega$ is flat.
\end {theorem}

Let us now consider a correspondence $\bb{X}$ 
defined over a number field $\mcal{F} \subset \bb{C}$. Let $\mathfrak{p}$ be a 
place of $\mcal{F}$ with $k_\mathfrak{p}$ as its residue field and 
$k_\mathfrak{p}^a$ its algebraic closure. Hence for each place $\mathfrak{p}$
we obtain a new correspondence $\mathbb{X}_\mathfrak{p} = (Y_\mathfrak{p}, 
X_\mathfrak{p}, \bar{\sigma}_{1,\mathfrak{p}},\bar{\sigma}_{2,\mathfrak{p}})
$ by base changing $\mathbb{X}$ to $k_\mathfrak{p}^a$. Also let us again
restrict ourselves to the case when $X=Y=\bb{P}^1$. Then by our Theorem 
\ref{mainthm2}, if for a place $\mathfrak{p}$, the group 
$\mcal{G}_{\bb{X}_\mathfrak{p}}$ is non-trivial, then the primitive $\omega$
of $\mcal{G}_{\bb{X}_\mathfrak{p}}$ has to be either $\omega = \frac{dt}{t-a}$,
the flat form of weight $1$ or $\omega = \frac{(dt)^2}{(t-a)(t-b)}$ which
is the flat form of weight $2$ provided that $deg~\sigma_1 \geq 
14 deg~\sigma_2$. Hence the above forms play a central role as
the semi-invariant forms for correspondences of $\bb{P}^1$ with unequal 
degrees.

Let us now look at an example. Consider $\bb{X} = (\bb{P}^1,\bb{P}^1,\sigma_1,
\sigma_2)$ defined over $\Z$. Let $\sigma_1 = t^m$ and $\sigma_2 = t^h$ and
$\omega_\mathfrak{p} = \frac{dt}{t}$ for all $\mathfrak{p}$. Then it is easy 
to see that for all $\mathfrak{p}$, $\omega_\mathfrak{p}$ is semi-invariant.
In fact, this simply follows because $\sigma_i^* \omega_\mathfrak{p} = \gamma_i
\omega_\mathfrak{p}$ for $i=1,2$ and some $\gamma_i$'s \cite{B1}. Now we
will construct another example by ``twisting'' the above example. Let $\sigma$
be any endomorphism of $\bb{P}^1$. Define $\sigma_1= \sigma \circ t^m$
and $\sigma_2 = \sigma \circ t^h$. For each $\mathfrak{p}$ define 
$\omega_\mathfrak{p}$ as before. Then clearly $\omega_\mathfrak{p}$ is 
semi-invariant because $\sigma_1^* \omega_\mathfrak{p} = \sigma^* (t^m)^* 
\omega_\mathfrak{p} = \sigma^* (\gamma_1 \omega_\mathfrak{p}) = 
\gamma_1 \sigma^* \omega_\mathfrak{p}$
and similarly $\sigma^*_2 \omega_\mathfrak{p} = \gamma_2 \omega_\mathfrak{p}$.

Our next result is the converse of the above example. It says that 
for a correspondence $\bb{X}=(\bb{P}^1, \bb{P}^1,\sigma_1,\sigma_2)$ defined 
over a number field $\mcal{F}$ and let $\bb{X}_\mathfrak{p}$ be the 
correspondence obtained from $\bb{X}$ as described above, if the group 
$\mcal{G}_\mathfrak{p}$ is non-trivial and is generated by primitive of weight 
$1$ for infinitely many places $\mathfrak{p}$, then that determines the 
morphisms $\sigma_1$ and $\sigma_2$ to be the ones described in the above 
example.
\begin {theorem}
\label{thm3}
Let $\bb{X}= (\bb{P}^1,\bb{P}^1,\sigma_1, \sigma_2)$ be defined over a 
number field $\mathcal{F}$ such that $\mathcal{G}_{\bb{X}_\mathfrak{p}}$
is non-trivial and is generated by primitives of weight $1$ for infinitely many
places $\mathfrak{p}$. Then there exists an endomorphism $\sigma:\bb{P}^1 \map 
\bb{P}^1$ defined over $\mathcal{F}$ such that,
\begin {equation}
(\sigma_1, \sigma_2) \sim (\lambda_1(\sigma \circ t^m), \lambda_2(\sigma \circ
t^h))  
\end {equation}
for some integers $m$ and $h$ and constants $\lambda_1$ and 
$\lambda_2$ as in corollary \ref{prodcor}.
\end {theorem}
We would now like to summarize our results and put them in context. The 
question that we are trying to solve was posed in \cite{B1}- classify all
$(\sigma_1,\sigma_2)$ of $\bb{X} = (\bb{P}^1,\bb{P}^1,\sigma_1,\sigma_2)$
defined over a number field $\mcal{F}$ such that $\mcal{G}_\bb{X_\mathfrak{p}}$
is non-trivial for infinitely many places $\mathfrak{p}$. The case when
$\sigma_1$ is fixed to be the identity morphism on $\bb{P}^1$ has been 
completely classified in \cite{B1}. Then it has been shown that the other 
morphism $\sigma_2 := \sigma$ can only be one of the {\it flat} maps- the 
{\it Multiplicative}, the {\it Chebyshev} and the {\it Latt\`{e}} maps. We 
will briefly describe the Multiplicative and Chebyshev maps for our purpose. A 
map $\sigma : \bb{P}^1 \map \bb{P}^1$ is multiplicative if $\sigma (t) =
t^{\pm d}$ for a positive integer $d$. A Chebyshev polynomial $\sigma$ of 
degree $d$ is the unique polynomial such that $\sigma(t+t^{-1}) =
t^d + t^{-d}$. Then for $\bb{X} = (\bb{P}^1,\bb{P}^1, \mathbbm{1},\sigma)$
which satisfies that $\mcal{G}_{\bb{X}_\mathfrak{p}}$ is non trivial for
infinitely many places $\mathfrak{p}$ then 
\begin {enumerate}
\item when the primitive of $\mcal{G}_{\bb{X}_\mathfrak{p}}$ is a flat form of 
weight $1$, then $\sigma$ is conjugate to a Multiplicative function.

\item when the primitive of $\mcal{G}_{\bb{X}_\mathfrak{p}}$ is a flat form
of weight $2$, then $\sigma$ is conjugate to a Chebyshev function.
\end {enumerate}
There is another possibility of the primitive other than the above two, 
as shown in \cite{B1}, which is associated with the Latt\`{e}s map. All these
above three described primitives exhausts all the possibilities for the group
$\mcal{G}_{\bb{X}_\mathfrak{p}}$ to be non-trivial for infinitely many 
$\mathfrak{p}$'s where $\bb{X}= (\bb{P}^1,\bb{P}^1,\mathbbm{1},\sigma)$.

Now in the general case where $\bb{X} = (\bb{P}^1,\bb{P}^1, \sigma_1,\sigma_2)$
is defined over $\mcal{F}$, if we insist that the pair of maps 
$(\sigma_1, \sigma_2)$ is conjugate to the ones totally ramified at one point,
say $\infty$, then our Theorem \ref{mainthm2} `morally' says that the 
primitive for a non-trivial $\mcal{G}_{\bb{X}_\mathfrak{p}}$ has to be a flat
form of weight either $1$ or $2$. However, we do indeed assume that extra
condition of $deg~\sigma_1 \geq 14 deg~\sigma_2$ to prove the result. One hopes
that this condition may be done away with or else there might be interesting
examples of primitives which are not flat forms in the case when $deg~\sigma_1
< 14 deg~\sigma_2$.

Our Theorem \ref{thm3} shows that if the primitive of 
$\mcal{G}_{\bb{X}_\mathfrak{p}}$ is a flat form of weight $1$, then the pair
of maps $(\sigma_1,\sigma_2)$ has to come from a pair of Multiplicative
functions composed with any arbitrary endomorphism $\sigma$ of $\bb{P}^1$.
Hence from the above analogy listed in the case of $(\mathbbm{1},\sigma)$, we
would like to conjecture that in the case when the primitive of 
$\mcal{G}_{\bb{X}_\mathfrak{p}}$ is a flat form of weight $2$, then our 
pair of morphisms $(\sigma_1,\sigma_2)$ come from a pair of Chebyshev functions
in the same manner as in the case of weight $1$. This question lies as one
of the motivation for future work for the author. Also the question of the 
third possibility of the primitive of $\mcal{G}_{\bb{X}_\mathfrak{p}}$ remains
completely open for further understanding.

{\bf Acknowledgements.} I would like to thank my advisor Prof. A. Buium for 
introducing me to the problem and for many helpful discussions. 
	
\section{\bf Proof of Theorem \ref{thm1}} 

\begin {proposition}
\label{cyclic}
If $deg~\sigma_1 > deg~\sigma_2$ then $\mcal{G}_{\bb{X}}$ is a free group 
of order $\leq 1$.
\end {proposition}
{\it Proof of proposition \ref{cyclic}.}
If $\mathcal{G}_{\bb{X}}= \{[1]\}$ then there is nothing more to prove. Hence
let us assume $\mathcal{G}_{\bb{X}}$ is non-trivial. 
Suppose $[\omega],[\omega'] \in \mathcal{G}_{\bb{X}}$ of the same weight $\nu$.
Then $f:= \frac{\omega}{\omega'} \in K(Y)$ is a rational function for $Y$.
But since we have $\sigma_1^* \omega = \lambda \sigma_2^* \omega$ and
$\sigma_1^* \omega' = \lambda' \sigma_2^* \omega'$ for some constants
$\lambda$ and $\lambda'$, we have
\begin {equation}
\sigma_1^*\left(\frac{\omega}{\omega'}\right) = \frac{\lambda}{\lambda'}
\sigma_2^* \left(\frac{\omega} {\omega'}\right)~ \imp ~
\sigma_1^* f = (\mbox{constant}) \sigma_2^* f
\end {equation} 
But $deg~\sigma_1^* f > deg~\sigma_2^* f$ unless $f$ is a constant, which 
implies that $[\omega] = [\omega']$. This shows that for a given weight $\nu$, 
there exists a unique class $[\omega] \in \mathcal{G}_{\bb{X}}$. 

Since $\mathcal{G}_{\bb{X}}$ is non-trivial, there exists a form $\omega$ 
with the smallest positive weight $\mu$ such that
$[\omega] \in \mathcal{G}_{\bb{X}}$.  And as we have shown above, this class of
weight $\mu$ is unique. Let $\omega'$ is a semi-invariant form of weight $\nu$. 
then $\nu = \mu l + r$ for some integer $r < \mu $. But then 
$\frac{\omega'}{\omega}$ is a semi-invariant form of weight $r$ which is a
contradiction to our hypothesis for $\omega$ unless $r=0~ \imp ~ [\omega'] = 
[\omega]^l$ and we are done. $\qed$

Let $\mathbb{X} = (Y,X,\sigma_1,\sigma_2)$ be a correspondence of
curves $X$ and $Y$ over any algebraically closed field $k$. Also  $\sigma_1$ 
and $\sigma_2$ are both {\it tamely ramified}. 
\begin {lemma}
\label {condlem}
If $\omega$ is an invariant form of weight $\nu$ of a correspondence 
$\mathbb{X}$ with $deg~\sigma_1 > deg~\sigma_2$ then,
\begin {equation}
\label{condbd}
\mathcal{C}_{Y,\omega} \leq \frac{1}{deg~\sigma_1 - deg~\sigma_2} 
(2deg~R_{\sigma_1} + deg~R_{\sigma_2})
\end {equation} 
where $R_{\sigma_1}$ and $R_{\sigma_2}$ are the ramification divisors of
$\sigma_1$ and $\sigma_2$ respectively.
\end {lemma}
\begin {remark}
Before proving the above lemma, we would like to remark that no assumption
on the smoothness of $X$ and $Y$ are necessary. Neither do we need to assume
any {\it properness} condition on either $X$ and $Y$. 
\end {remark}

{\it Proof of lemma $\ref{condlem}$.} We denote, $d_i = deg~\sigma_i$ for 
$i=1,2$.
Now let us write the divisor associated to our invariant form $\omega$ as,
$$ div~(\omega) = \sum_{i=1}^{\mathcal{C}_{Y,\omega}} f_iy_i, f_i \in 
\mathbb{Z} \backslash (0),~ y_i \in Y$$
where $\mathcal{C_{Y,\omega}}$ is the conductor of a form $\omega$ as defined 
before. Then for each $y_i$ define its pull-back via $\sigma_1$ as
\beqar
\sigma^*_1y_i & = & \sum_{j=1}^{n_i} e_{ij} \beta_{ij},~ \mbox{such that}\\
& & e_{ij} \leq 1,\mbox{ when }  j \leq \nu_i,~\mbox{for some}~ \nu_i 
\mbox{ and}\\
& & e_{ij} \geq 2,\mbox{ when }  j > \nu_i
\eeqar 
In particular, for all $j \leq \nu_i$, the points $\beta_{ij} \in X$ are
unramified under the map $\sigma_1$. Then one can write the ramification 
divisor for $R_{\sigma_1}$ as
\begin {equation}
 \label{ramdiv}
 R_{\sigma_1} = \sum_{i=1}^{\mathcal{C}_{Y,\omega}} \sum_{j=1}^{n_i} (e_{ij}-1)
\beta_{ij} + \sum_k l_k \delta_k
\end {equation}
where $\delta_k$'s are the other ramification points with ramification indices
$l_k$ in $X$ which do not belong to the preimages of $y_i$'s. Hence writing
out $div~(\sigma_1^*\omega)= \sigma_1^*div~(\omega)+ \nu R_{\sigma_1}$,
we obtain,
\beqar
div~(\sigma_1^*\omega) & = & \sum_{i=1}^{\mathcal{C}_{Y,\omega}} 
\sum_{j=1}^{n_i} (f_i e_{ij} + \nu(e_{ij}-1))\beta_{ij} + \sum_k l_k \delta_k \\
& = & \sum_{i=1}^{\mathcal{C}_{Y,\omega}}\left( \sum_{j=1}^{n_i}(e_{ij}(f_i+\nu)
- \nu) \right) \beta_{ij} + \sum_k l_k \delta_k  
\eeqar
Note that if $j \leq \nu_i$, then $\beta_{ij} \in \mbox{support }
(div~(\sigma_1^* \omega)) = \mbox{support }(div~(\sigma_2^*\omega))$. 
Hence we have
\begin {equation}
\label {supp}
 \sum_{i=1}^{\mathcal{C}_{Y,\omega}} \nu_i \leq |\mbox{support }
(div~(\sigma_2^*\omega))|
\end {equation}
We have,
\beqar
\mathcal{C}_{Y,\omega}d_1 & = & \sum_{i=1}^{\mathcal{C}_{Y,\omega}} 
\sum_{j=1}^{n_i}
e_{ij} \\
& = & \sum_{i=1}^{\mathcal{C}_{Y,\omega}} \left( \sum_{j=1}^{\nu_i} 1 + 
\sum_{j=\nu_i+1}^{n_i} e_{ij} \right) \\
& = & \sum_{i=1}^{\mathcal{C}_{Y,\omega}} \nu_i + 
\sum_{i=1}^{\mathcal{C}_{Y,\omega}}\sum_{j=\nu_i+1}^{n_i} e_{ij}  
\eeqar
and hence,
\begin {equation}
\label{cond1}
\mathcal{C}_{Y,\omega}d_1 \leq  |\mbox{support } (div~(\sigma_2^*\omega))| + 
\sum_{i=1}^{\mathcal{C}_{Y,\omega}}\sum_{j=\nu_i+1}^{n_i} e_{ij} 
\end {equation}
From \ref{ramdiv} we have, 
\beqar
deg~R_{\sigma_1} & = & \sum_{i=1}^{\mathcal{C}_{Y,\omega}} \sum_{j=\nu_i}^{n_i}
(e_{ij} - 1) + \sum_k l_k \\
& \geq & \sum_{i=1}^{\mathcal{C}_{Y,\omega}} \sum_{j=\nu_i}^{n_i} e_{ij} - 
\sum_{i=1}^{\mathcal{C}_{Y,\omega}} \sum_{j=\nu_i}^{n_i} 1 \\
deg~R_{\sigma_1} + \sum_{i=1}^{\mathcal{C}_{Y,\omega}} \sum_{j=\nu_i}^{n_i} 1 
& \geq & \sum_{i=1}^{\mathcal{C}_{Y,\omega}} \sum_{j=\nu_i}^{n_i} e_{ij} \\
\eeqar
\begin {equation}
\label{cond2}
2 deg~R_{\sigma_1}  \geq  \sum_{i=1}^{\mathcal{C}_{Y,\omega}} \sum_{j=\nu_i}^{n_i} e_{ij},~\mbox{because } \sum_{i=1}^{\mathcal{C}_{Y,\omega}} \sum_{j=\nu_i}^{n_i} 1 \leq deg~R_{\sigma_1}    
\end {equation}
Now $|\mbox{support } (div~(\sigma_2^*\omega))| \leq |\mbox{support } 
(\sigma_2^*(div~(\omega)))| + |\mbox{support }\nu(R_{\sigma_2})|$ which gives,
\begin {equation}
\label{cond3}
|\mbox{support } (div~(\sigma_2^*\omega))| \leq \mathcal{C}_{Y,\omega} d_2 +
 deg~ R_{\sigma_2}
 \end {equation}  

Therefore, putting \ref{cond1}, \ref{cond2} and \ref{cond3} together we obtain,
\begin {equation}
\mathcal{C}_{Y,\omega} \leq \frac{1}{d_1-d_2} (2deg~R_{\sigma_1} + 
deg~R_{\sigma_2})
\end {equation}
and we are done. $\qed$ 

As a corollary, we prove our first theorem,

{\it Proof of theorem 1.}
Since our curves $X$ and $Y$ are smooth and projective with $\sigma_1$ and 
$\sigma_2$ tamely ramified, by Riemann-Hurwitz \cite{Hartshorne}  we get
\begin {equation}
\label{deggen}
deg~R_{\sigma_i} = (2g_X-2) - d_i(2g_Y-2),~\forall i = 1,2
\end {equation}
Hence combining lemma \ref{condlem} and \ref{deggen} we obtain our result.
$\qed$

\section{\bf Proof of Theorem \ref{mainthm2}} 
We consider the restriction of our maps $\sigma_1,~\sigma_2$ and $\omega$ to
the affine line $\bb{A}^1 \subset \bb{P}^1$ to obtain a new correspondence
$\bb{X} = (\bb{A}^1, \bb{A}^1,\sigma_1,\sigma_2)$. We will also assume that 
$deg~\sigma_1 > deg~\sigma_2$. 
  
Let $div~(\omega) = \sum_{i=1}^m e_i y_i$. Then we have the following lemma,
\begin {lemma}
If $\mathbb{X}$ admits a semi-invariant form $\omega$ of weight $\nu$ then
$$\sum_{i=1}^{\mathcal{C}_{\bb{A}^1,\omega}} e_i = -\nu$$
\label {lem1}
\end {lemma}
{\it Proof.} Since $\sigma_i$'s are endomorphisms of $\mathbb{A}^1$, we 
  note that the degree of the ramification divisor 
$R_{\sigma_i} = d_i-1$ for $i=1,2$. Hence equating the degrees of the 
divisors $div~(\sigma^*_1 \omega)$ and $div~(\sigma^*_1 \omega)$ we obtain,
\beqar
\sum_{i=1}^{\mathcal{C}_{\bb{A}^1,\omega}} e_i d_1 + \nu(d_1-1) & = & 
\sum_{i=1}^{\mathcal{C}_{\bb{A}^1,\omega}} e_i d_2 + \nu(d_2-1),\\
\sum_{i=1}^{\mathcal{C}_{\bb{A}^1,\omega}} e_id_1 + \nu d_1 - \nu & = & 
\sum_{i=1}^{\mathcal{C}_{\bb{A}^1,\omega}} e_id_2 + \nu d_2 - \nu,\\
(d_1 - d_2) \sum_{i=1}^{\mathcal{C}_{\bb{A}^1,\omega}} e_i & = & \nu 
(d_2 - d_1), \\
\sum_{i=1}^{\mathcal{C}_{\bb{A}^1,\omega}} e_i & = & - \nu ~~ \qed
\eeqar

\begin {lemma}
If $\mathbb{X}$ admits a non-zero semi-invariant form $\omega$ then
$$ \mathcal{C}_{\bb{A}^1,\omega} \leq 2+ \frac{3(d_2-1)}{d_1-d_2}$$
\end {lemma}
{\it Proof.} Since the degree of the ramification 
divisors $R_{\sigma_1}$ and $R_{\sigma_2}$ are $d_1-1$ and $d_2-1$ respectively, substituting in \ref{condbd} we obtain our desired result. $\qed$    

\begin {corollary}
\label {cor1}
If $d_1 \geq 4d_2$ then $\mathcal{C}_{\bb{A}^1,\omega} \leq 2$.
\end {corollary}
{\it Proof.}
\beqar
d_1 & \geq & 4d_2 \\
d_1 - d_2 & \geq & 3d_2 \\
\frac {1}{d_1-d_2} & \leq & \frac{1}{3d_2} \\
\frac {3(d_2-1)}{d_1-d_2} & \leq & \frac{3d_2-3}{3d_2} \\
& = & 1 - \frac{1}{d_2} \\
& < & 1 \\
2+ \frac{3(d_2-1)}{d_1-d_2} & < & 3
\eeqar
In other words $\mathcal{C}_{\bb{A}^1,\omega}< 3$ but since 
$\mathcal{C}_{\bb{A}^1,\omega}$ is a natural number we conclude
that $\mathcal{C}_{\bb{A}^1,\omega} \leq 2$. $\qed$ 

\begin {lemma}
\label{1case}
Let $\mathbb{X}$ be a correspondence such that $d_1 \geq 4d_2$ and admitting
a semi-invariant form  $\omega$ with $\mathcal{C}_{\bb{A}^1,\omega}=1$, then 
$\omega$ is a flat form of type-$1$.
\end {lemma}

{\it Proof.} Since $\mathcal{C}_{\bb{A}^1,\omega}=1$ and by lemma \ref{lem1}, 
$div~(\omega) = -\nu .y$ which is precisely a flat form of type-$1$. $\qed$   

\begin {lemma}
\label {weight1}
If $\omega$ is a semi-invariant form with $\mathcal{C}_{\bb{A}^1,\omega}=2$ 
then its weight cannot be $1$.
\end {lemma}
{\it Proof.} By lemma \ref{lem1} we may assume that 
\begin {equation}
div~\omega = ey_1 -(e+1)y_2,\mbox{ and } e \geq 1. 
\end {equation}
Let the pull-pack of the divisor $y_2$ via the two maps be,
\begin {equation}
\sigma_1^*y_2 = \sum_{i=1}^n e_i \beta_i \mbox{ and } 
\sigma_2^*y_2 = \sum_{j=1}^m f_j \beta_j' 
\end {equation} where $e_i$ and $f_j$'s are positive non-zero integers. 
Then one can write the ramification divisors as in (\ref{ramdiv}) as,
\begin {equation}
R_{\sigma_1} = \sum_{i=1}^n (e_i-1) \beta_i + D \mbox{ and }
R_{\sigma_2} = \sum_{j=1}^m (f_j-1) \beta_j' + D'
\end {equation}  for some effective divisors $D$ and $D'$. Hence, if we 
equate the poles of the divisor $div~\sigma_1^*\omega$ and $div~\sigma_2^*
\omega$ we obtain,
\beqar
\sum_{i=1}^n ((e_i-1)-e_i(e+1)) \beta_i & = & \sum_{j=1}^m ((f_j-1)-f_j(e+1)) 
\beta_j' \\
\sum_{i=1}^n -(1+e_ie)\beta_i & = &  \sum_{j=1}^m -(1+f_if)\beta_j'
\eeqar  
Since $-(1+e_i e)$ and $-(1+f_j f)$ are non-zero for any $i$ and $j$
implies that $n = m,~ \beta_i = \beta_j'$ and $e_i = f_j$. But then 
$deg~\sigma_1 = deg~\sigma_2$ which is a contradiction. $\qed$ 

\begin {theorem}
\label{thm2}
If $d_1 \geq 14 d_2$ and $\omega$ is a semi-invariant form for 
$\mathbb {X}$ such that $[\omega]$ is the primitive and 
$\mathcal{C}_{\bb{A}^1,\omega} =2$, then $div~(\omega) = -y_1-y_2$. 
\end {theorem}

{\it Proof.} Let $div~(\omega) = e.y_1 + f.y_2$ for some $e,f \in \mathbb{Z}
\backslash (0)$. Then by lemma \ref{lem1}, $e+f = -\nu$ and $(e,f) =1$
because $\omega$ is primitive.  We also know that
$$ div~(\sigma_i^*\omega) = \sigma_i^*(div~(\omega)) + \nu.R_{\sigma_i}~
\mbox{, for } i = 1,2$$ 
Then we have,
$$\sigma_1^*(y_1) = \sum_{i=1}^n e_i\alpha_i ~~ \mbox{and } 
\sigma_2^*(y_1) = \sum_{i=1}^{n'} e_i'\alpha_i$$
$$\sigma_1^*(y_2) = \sum_{j=1}^m f_j\beta_j ~~ \mbox{and } 
\sigma_2^*(y_2) = \sum_{j=1}^{m'} f_j'\beta_j'$$
for some $e_i,e_i',f_j,f_j'> 0 $. We can also write the ramification
divisors as-
$$R_{\sigma_1} = \sum_{i=1}^n(e_i-1)\alpha_i + \sum_{j=1}^m (f_j-1) \beta_j
+ \sum_{k=1}^p l_k\delta_k$$
$$R_{\sigma_1} = \sum_{i=1}^{n'}(e_i'-1)\alpha_i' + \sum_{j=1}^{m'} (f_j'-1) 
\beta_j+ \sum_{k=1}^p l_k'\delta_k'$$ 
Then the divisor associated to $\sigma_1^*\omega$ is
\beqar
div~(\sigma_1^*\omega) & = & e\sum_{i=1}^n e_i\alpha_i + f\sum_{j=1}^m f_j
\beta_j + \nu\left(\sum_{i=1}^n(e_i-1)\alpha_i + \sum_{j=1}^m (f_j-1) \beta_j
+ \sum_{k=1}^p l_k\delta_k \right) \\
& = & \sum_{i=1}^n (ee_i+\nu e_i - \nu) \alpha_i + \sum_{j=1}^m (ff_j+ \nu f_j
-\nu)\beta_j + \sum_{k=1}^p \nu l_k \delta_k \\
& = & \sum_{i=1}^n (-fe_i-\nu) \alpha_i + \sum_{j=1}^m (-ef_j-\nu)
\beta_j + \sum_{k=1}^p \nu l_k \delta_k, ~~ \mbox{since } e+f = -\nu 
\eeqar
And similarly, we have
$$div~(\sigma_2^*\omega) = \sum_{i=1}^{n'}(-fe_i'-\nu) \alpha_i' + \sum_{j=1}
^{m'} (-ef_j'-\nu) \beta_j' + \sum_{k=1}^{p'} \nu l_k' \delta_k' $$

{\it Claim. If any one of the factors of the form $(-fe_i -\nu)$ or 
$(-fe_i'-\nu)$ is $0$ then $f =-1$.}

{\it Proof.} If $-fe_i-\nu = 0$ for some $i$ implies that $f=\frac{\nu}{-e_i}$.
That means that $f$ is negative because both $\nu$ and $e_i$ are positive.
By lemma \ref{lem1} we know that $f+e = f(-e_i)=\nu$ which implies that,
$e= -f(1-e_i)~ \imp~ f~|~e$. But since $\omega$ is primitive, we have $(e,f)
=1$. Hence $f$ can only be $1$ or $-1$ but then $f$ is negative 
and hence $f =-1$. Similar argument implies the result in the case when
$-fe_i'-\nu = 0$ and this completes the proof of our claim.

Similarly we can show that,

{\it Claim. If any one of the factors of the form $(-ef_j -\nu)$ or 
$(-ef_j'-\nu)$ is $0$ then $e =-1$.}

Now suppose $\omega$ is not a flat form. Then from 
 the above claims, we have the following two cases to consider:

\underline {Case (1):}

If none of the $e$ and $f$ equals $-1$. Then the above claim implies that
$ (-fe_i-\nu),~ (-fe_i'-\nu)  \ne 0,~ \forall~i $ and  $(-ef_j-\nu),~(ef_j'-\nu)
\ne 0,~\forall~j$ which means, all the coefficients in $div~(\sigma_1^* \omega)$
and $div~(\sigma_2^* \omega)$ are non-zero and since we have,
\begin {equation}
|\mbox{support }div~(\sigma_1^*\omega)| = |\mbox{support }div~(\sigma_2^*
\omega)|
\end {equation}
imples that $m+n+p = m'+n'+p'$.

\underline {Case (2):}
We may assume without loss of generality that $e=-1,~ \imp~ f = 1-\nu$ and 
$f \ne -1$ as because then our $\omega$ is already in the required form. Since
$f \ne -1$ then by the above claim, $(-fe_i-\nu)\ne 0,~ \forall~ i$ and 
$(-fe_i'-\nu) \ne 0,~\forall i$

If \underline {Case (1):} was true then from the fact,
$$ deg~(div~(\sigma^*_1\omega)) =  deg~(div~(\sigma^*_2\omega))$$
we obtain,
\beqar
\sum_{i=1}^n(-fe_i-\nu)+\sum_{j=1}^m(-ef_j-\nu)+\sum_{k=1}^p\nu l_k
& = & \sum_{i=1}^{n'}(-fe_i'-\nu)+\sum_{j=1}^{m'}(-ef_j'-\nu)+\sum_{k=1}^{p'}
\nu l_k' \\
-f\sum_{i=1}^n e_i-\sum_{i=1}^n \nu -e\sum_{i=1}^m f_j -\sum_{j=1}^m\nu +
\nu \sum_{k=1}^pl_k  & = &
-f\sum_{i=1}^{n'} e_i'-\sum_{i=1}^{n'} \nu -e\sum_{i=1}^{m'} f_j' -
\sum_{j=1}^{m'}\nu + \nu \sum_{k=1}^{p'}l_k' \\
-fd_1-\nu n - ed_1-\nu m + \nu\sum_{k=1}^pl_k & = &
-fd_2-\nu n' - ed_2-\nu m' + \nu\sum_{k=1}^{p'}l_k'\\
(-f-e)d_1 -\nu(n+m) +\nu \sum_{k=1}^p l_k & = &
(-f-e)d_2 -\nu(n'+m') +\nu \sum_{k=1}^{p'} l_k'  \\
\nu d_1 -\nu(n+m) + \nu\sum_{k=1}^p l_k & = &
\nu d_2 -\nu(n'+m') + \nu\sum_{k=1}^{p'}l_k' ~~~ \dots (*)
\eeqar
Now the right-hand side is bounded from above by
\beqar
\nu d_2-\nu(n'+m')+ \nu \sum_{k=1}^{p'} & \leq & \nu d_2 + \nu \sum_{k=1}^{p'}
l_k' ~~~~~~~~~~~ (\mbox{ Since }  \nu(n'+m') \geq 0 ) \\
& \leq & \nu d_2 +\nu (d_2-1) ~~~~~~~ (\mbox{ Since } \sum_{k=1}^{p'}
l_k' \leq d_2-1) \\
& = & \nu(2d_2-1) \hspace{2cm} \dots(1)
\eeqar
Next we can find a lower bound for the left-hand side in the above inequality 
as follows,
$$ \nu d_1 - \nu(n+m) +\nu\sum_{k=1}^pl_k~ \geq ~\nu d_1 -\nu(n+m)
~~~~~~~(\mbox{ Since }\nu\sum_{k=1}^p l_k \geq 0) $$
Now
\beqar
 n+m & \leq & m+n+p \\
 & = & m'+n'+p' \\
 & \leq & d_2+d_2+(d_2-1) ~~~~~~(\mbox{ Since } m'\leq d_2,~n'\leq d_2,~
p' \leq d_2-1~) \\
 & = & 3d_2-1
\eeqar
And hence we get
\beqar
\nu d_1 -\nu(n+m) +\nu \sum_{k=1}^p l_k & \geq & \nu d_1-\nu(3d_2-1)\\
& = & \nu(d_1-3d_2-1)\\
& \geq & \nu(14d_2-3d_2 -1)~~~~~~~~~(\mbox{ Since } d_1 \geq 14d_2~)\\
& = & \nu(11d_2-1) \hspace{2cm} (2)
\eeqar
Therefore combining (1) and (2) we get
$$ \nu d_1 - \nu(n+m) + \nu\sum_{k=1}^pl_k~ \geq~ \nu(11d_2-1) ~ \gneqq ~
\nu(2d_2-1) \geq \nu d_2 -\nu(n'+m') +\nu\sum_{k=1}^{p'}l_k'$$
which contradicts $(*)$ and hence removes the possibility of
\underline{Case (1)} for $\omega$.

Now we proceed to show that \underline{ Case (2)} is also not possible for
$\omega$. In this situation as discussed above, we have $e=-1$ and $f = 1-\nu 
\ne -1$.

For any divisor $H = \sum a_i P_i$, we define $|H| = \sum |a_i|$. Then 
consider the divisor
$$D = \sum_{i=1}^n (-fe_i'-\nu)\alpha_i + \sum_{k=1}^{p'}\nu l_k\delta_k$$
Then we have, $deg~D \leq |D| \leq |div~\sigma_1^*\omega| = |div~\sigma_2^*
\omega|$.  
\beqar
deg~D & \leq & \sum_{i=1}^{n'}|(-fe_i'-\nu)|+\sum_{j=1}^{m'}|(-ef_j'-\nu)|
+\sum_{k=1}^{p'}|\nu l_k'|\\
& \leq & \sum_{i=1}^{n'}(|-fe_i'|+|\nu |) +\sum_{j=1}^{m'}(|-ef_j'|+|\nu|)
+\sum_{k=1}^{p'}\nu l_k' ~~~~~~~ (\mbox{ Since } l_k' \geq 0,~\forall~k ~)
\eeqar
Since $f < 0$ we get that $-fe_i' > 0,~ \forall~i$ and hence $|-fe_i'|=-fe_i'$.
and also by the same reason we obtain $|-ef_i'|=-ef_j',~\forall~j$. And
$\nu > 0$, our above inequality becomes:
\beqar
deg~D & \leq & \sum_{i=1}^{n'}(-fe_i'+\nu ) +\sum_{j=1}^{m'}(-ef_j'+\nu)
+\sum_{k=1}^{p'}\nu l_k' ~~~~~~~ (\mbox{ Since } l_k' \geq 0,~\forall~k ~)\\
& = & -f\sum_{i=1}^{n'}e_i' + \sum_{i=1}^{n'}\nu -e\sum_{j=1}^{m'}f_j'
+\sum_{j=1}^{m'}\nu +\nu \sum_{k=1}^{p'}l_k'\\
&=& -fd_2+\nu(n'+m') -ed_2+\nu\sum_{k=1}^{p'}l_{k'}\\
& =& (-e-f)d_2+\nu(n'+m')+\nu\sum_{k=1}^{p'}l_k'\\
& \leq & \nu d_2 +2d_2\nu +\nu (d_2-1) ~~~~(\mbox{Since }n',~m' \leq d_2
\mbox{ and } \sum_{k=1}^{p'}l_k' \leq d_2-1~) \\
&=& \nu (d_2+2d_2+d_2-1) \\
& = & \nu(4d_2-1) ~~~~~~~~(3)
\eeqar
Also, considering the degree of the divisor $D$ we obtain,
\beqar
deg~D & = & \sum_{i=1}^n (-fe_i-\nu) +\sum_{k=1}^p\nu l_k\\
& = &-f\sum_{i=1}^n e_i -\sum_{i=1}^n \nu + \nu\sum_{k=1}^pl_k \\
& \geq & -fd_1 -\nu n ~~~ (\mbox{ Since } \sum_{k=1}^p l_k\geq 0~)\\
& = & (\nu-1)d_1-\nu n
\eeqar
Here we note that $n \leq n'+m'+p'$ since none of the $(-fe_i-\nu)$'s are 
zero and the right hand side of the inequality is the size of the support of
the divisor $div~(\sigma_2^*\omega)$ and we know $n',~m' \leq d_2$
and $p' \leq d_2-1$. Hence we get that $n \leq 3d_2-1$. Substituting this
in the above inequality we obtain,
\beqar
deg~D &\geq & (\nu-1)d_1 -\nu(3d_2-1) \\
& = & \nu((1-\frac{1}{\nu})d_1 - 3d_2+1)
\eeqar
By lemma \ref{weight1}, $\nu \geq 2 ~\imp~ \frac{1}{\nu}\leq\frac{1}{2}~\imp~
-\frac{1}{\nu} \geq-\frac{1}{2}~ \imp~ 1-\frac{1}{\nu}\geq 1-\frac{1}{2}=
\frac{1}{2}$. Hence we have
$$ deg~D \geq \nu(\frac{1}{2}d_1 -3d_2+1)$$
But by our hypothesis, $d_1 \geq 14d_2$ and we get,
$$deg~D \geq \nu(\frac{1}{2}.14d_2-3d_2+1) = \nu(4d_2+1)~~~~~~~(4)$$
We combine (3) and (4) to obtain
$$deg~D ~ \geq~\nu(4d_2+1) \gneqq \nu(4d_2-1) \geq deg~D$$
and here lies the contradiction for \underline{Case (2)}. $\qed$

{\it Proof of theorem \ref{mainthm2}.} 
By corollary \ref{cor1} we need to check for the two cases when 
$\mathcal{C}_{\bb{A}^1,\omega}$ is either $1$ or $2$.  
Lemma \ref{1case} shows that when $\mathcal{C}_{\bb{A}^1,\omega} =1$ then
$\omega = \frac{1}{t-y}dt$, which implies that weight of $\omega$ is $1$ and
$div~\omega = -y-\infty$.    

And when $\mathcal{C}_{\bb{A}^1,\omega} =2$, theorem \ref{thm2} shows that 
$\omega = \frac{1}{(t-y_1)(t-y_2)}(dt)^2$, which implies that the weight of
$\omega$ is $2$ and $div~\omega= -y_1-y_2$. And this ends the proof. 
$\qed$  

\section {\bf Proof of theorem \ref{thm3}}
Call $\tilde{\omega} = \sigma_1^* \omega = \sigma_2^* \omega$. Let $x \in X$ 
then by \cite{B1},
\begin {equation}
\label{ordlemma}
ord_x(\tilde{\omega}) + \nu = e_{\sigma_i}(x)(ord_{\sigma_i(x)}\omega +\nu)
~\mbox{for } i = 1,2
\end {equation}

If $\omega$ is a flat form of weight $1$, it is easy to see that there exists
an automorphism of $\bb{P}^1$ fixing $\infty$ such that $\sigma_1' \sim 
\sigma_1$ and $\sigma_2' \sim \sigma_2$ and $\frac{dt}{t}$ is a 
semi-invariant form for $(\sigma_1',\sigma_2')$. Hence it is sufficient to 
assume that the semi-invariant form $\omega$ is of the form $\frac{dt}{t}$. 

\begin {proposition}
\label {eqsupp}
When $\nu = 1$, then $\mbox{support }\sigma_1^*(div~\omega) = \mbox{support }
\sigma_2^*(div~\omega)$.
\end {proposition}
{\it Proof.} 
Let $\omega$ be the primitive semi-invariant form associated to 
$\mathcal{C}_{\bb{A}^1,\omega} = 1$. Then $div~(\omega) = -0 - \infty \in 
div~(\bb{P}^1)$. Then,
\beqar
y \in \mbox{support }\sigma_1^*(div~\omega) & \iff & \sigma_1(y)  \in 
\mbox{support }div~\omega \\
\imp ord_{\sigma_1(y)}\omega & = & -1 \\
\imp e_{\sigma_2(y)}(ord_{\sigma_2(y)}\omega+1) &=& 0, \mbox{ because of 
\ref{ordlemma} }\\
\imp ord_{\sigma_2(y)}\omega & =& -1 \mbox{ since }e_{\sigma_2}(y) \geq 1
\eeqar
which means that $y \in \mbox{support }\sigma_2^*(div~\omega)~\imp~
\mbox{support }\sigma_1^*(div~\omega) \subset \mbox{support } \sigma_2^*
(div~\omega)$ and similarly we can show for the other direction and that
concludes our proof. $\qed$

Let $S = \{\alpha_1,...,\alpha_n\} = \mbox{support }\sigma_i^*\{0\}$ 
for $i=1,2$.
\begin {corollary}
\label {prodcor}
We have the following expressions for our $\sigma_1$ and $\sigma_2$,
\begin {equation}
 \sigma_1(t) = \lambda_1\prod_{i=1}^n (t- \alpha_i)^{e_i}\mbox{ and } 
 \sigma_2(t) = \lambda_2\prod_{i=1}^n (t- \alpha_i)^{f_i}
\end {equation}
 for some tuple of positive non-zero integers $(e_1,...,e_n)$ and 
$(f_1,...,f_n)$ and $\lambda_1$ and $\lambda_2$ are constants in the 
ground ring.
\end {corollary}
{\it Proof.} This follows from proposition \ref{eqsupp}. Since both 
$\sigma_1$ and $\sigma_2$ are totally ramified at $\infty$, all the elements 
of $S$ are all the roots of the polynomials describing the maps $\sigma_1$
and $\sigma_2$ respectively. Hence we obtain our required factorisation 
of both $\sigma_1$ and $\sigma_2$ and we are done. $\qed$


Let $k$ be a field of characteristic $p$. Then there exists a group homomorphism
$\theta_p: \bb{Z}_{(p)}^\times \map k^\times$ defined by 
$\theta_p \left(\frac{m}{n}\right) = 
\frac{\bar{m}}{\bar{n}}$ where $\bar{m}$ denotes the reduction of $m$ modulo $p$.

\begin {lemma}
\label{ker}
$ker~ \theta_p = \left\{\frac{m}{n}~|~ m \equiv n~\mbox{ mod } p\right\}$
\end {lemma}
{\it Proof.} $\frac{m}{n} \in ker~\theta_p$ iff $\frac{\bar{m}}{\bar{n}} = 1$ iff
$\bar{m} \equiv \bar{n}$ and we are done. $\qed$

Let $\bb{X}$ now be a correspondence defined over a number field $\mathcal{F}$.
Let $p$ be the characteristic of the residue field at $\mathfrak{p}$.

\begin {lemma}
If $(2deg~\sigma_1)(deg~\sigma_2) < p$ then there exists a function $\sigma$ 
defined over $\mathcal{F}$ such that,
$$\sigma_1 (t)= \lambda_1(\sigma(t))^m \mbox{ and } \sigma_2(t) = 
\lambda_2(\sigma(t))^h$$ for some non-negative
integers $m$ and $h$ and $\lambda_1$ and $\lambda_2$ as in corollary 
\ref{prodcor}. In other words, both $\sigma_1$ and $\sigma_2$ are
a composition of $\sigma$ and a multiplicative function $\lambda_1t^m$ or 
$\lambda_2 t^h$ respectively. 
\end {lemma}

{\it Proof.} 
\beqar
\sigma^*_1(\omega) & = &l \sigma^*_2(\omega), \mbox{ for some }l\\
\frac{\sigma_1'}{\sigma_1} & = & l\frac{\sigma_2'}{\sigma_2} \\
\sum_{i=1}^n \frac{e_i}{t-\alpha_i} & = & l \sum_{i=1}^n
\frac{f_i}{t-\alpha_i} \mbox{ by corollary } \ref{prodcor}
\eeqar
Note that $(e_i,p)= (f_i,p)=1$ for all $i$ because $e_i <d_1<p$ and $f_i<d_2<p$.
Hence the above equality is possible only if $\theta_p(\frac{e_i}{f_i}) = 
\theta_p 
(\frac{e_j}{f_j}) \equiv l$ for all $i$ and $j$, i.e $ \frac{e_if_j}{f_ie_j} \in
ker~\theta_p$. By lemma \ref{ker} we have $e_if_j - f_je_i = pg$ for some integer
$g$. We claim that $g = 0$. If not then $|e_if_j-f_ie_j| =|p| |g| > |p|$. On the
other hand, $|e_if_j-f_ie_j| \leq |e_if_j| + |e_jf_i| \leq 2d_1 d_2$ but 
$2d_1d_2 < p$ which is a contradiction and proves our claim. In other words we
have $\frac{e_i}{f_i} = \frac{e_j}{f_j} = \frac{m}{h}$ for all $i$ and $j$
and for some integers $m$ and $h$ with $(m,h) = 1$. Hence we have $e_ih = f_im$
 for all $i$. Take $g_i = e_i/m = f_i/h$ and set 
$\sigma(t) = \prod_{i=1}^n (t-\alpha_i)^{g_i}$.  Then
one checks easily that $\sigma_1(t) = \lambda_1(\sigma(t))^m$ and $\sigma_2(t) 
= \lambda_2(\sigma (t))^h$ and we are done. $\qed$

{\it Proof of theorem \ref{thm3}.}
Since $\bb{X}$ over the number field $\mathcal{F}$ admits the flat form
for infinitely many places, chose a $\mathfrak{p}$ such that its residue
field has characteristic $p > 2 deg~\sigma_1 deg~\sigma_2$. Then the $e_i$'s and
$f_i \in \bb{Z}_{(p)}^\times$ and the theorem follows. $\qed$    

\bibliographystyle{amsplain}

\vspace {1cm}

{\tiny UNIVERSITY OF NEW MEXICO, ALBUQUERQUE, NM 87131}

{\it \small E-mail address:} {\small \tt arnab@math.unm.edu}
\end {document}